\documentclass[12pt]{article}

\usepackage{amsmath}
\usepackage{latexsym}
\usepackage{amssymb,enumerate}
\usepackage{amsfonts,hyperref}
 
\usepackage{rotating}
\usepackage{graphics}
\include{PDF}

\newcommand{\proof}{\par\noindent{\it Proof.\ \ }}
\def\qed{\ifmmode\square\else\nolinebreak\hfill
$\Box$\fi\par\vskip12pt}

\DeclareMathOperator\notdiv{{\,\not\big|\,}}
 
\def\calD{{\mathcal D}} 
\def\calP{{\mathcal P}}
\def\calS{{\mathcal S}} 
\def\calB{{\mathcal B}} 
\def\calM{{\mathcal M}} 
\def\FF{\mathbb F} 
\def\ZZ{\mathbb Z} 
\def\HS{{\rm HS}} 
\def\GF{{\sf GF}}

\DeclareMathOperator\PG{{\sf PG}}
\DeclareMathOperator\AG{{\sf AG}}

\DeclareMathOperator\soc{{\sf soc}}

\DeclareMathOperator\Aut{{\sf Aut}}

\DeclareMathOperator\Sp{{\sf Sp}}
\DeclareMathOperator\GammaL{{\sf \Gamma L}}
\DeclareMathOperator\AGammaL{{\sf A\Gamma L}}
\DeclareMathOperator\PGammaL{{\sf P\Gamma L}} \DeclareMathOperator\PSigmaL{{\sf P\Sigma L}}
\DeclareMathOperator\PSL{{\sf PSL}}\DeclareMathOperator\PGL{{\sf PGL}}
\DeclareMathOperator\GL{{\sf GL}} \DeclareMathOperator\SL{{\sf SL}}
\DeclareMathOperator\AGL{{\sf AGL}} 
\DeclareMathOperator\Sym{{\sf Sym}}\DeclareMathOperator\Alt{{\sf Alt}} 
\DeclareMathOperator\ASL{{\sf ASL}}

 \DeclareMathOperator\PSU{{\sf PSU}}   \DeclareMathOperator\PGammaU{{\sf P\Gamma U}} 
 \DeclareMathOperator\G{{\sf G}}
\DeclareMathOperator\Ree{{\sf Ree}}

\newtheorem{theorem}{Theorem}[section]%
\newtheorem{lemma}[theorem]{Lemma}%
\newtheorem{corollary}[theorem]{Corollary}%

\newtheorem{proposition}[theorem]{Proposition}%
\newtheorem{definition}[theorem]{Definition}%
\newtheorem{example}[theorem]{Example}%
\newtheorem{construction}[theorem]{Construction}%
\newtheorem{remark}[theorem]{Remark}%

\input amssym.def
\input amssym.tex

\begin{document}

\title{Pairwise transitive 2-designs}
\author{Alice Devillers and Cheryl E.~Praeger\thanks{This  research forms part of the Discovery Project DP130100106 of the authors funded by the Australian Research Council.} \\ Centre for the Mathematics of Symmetry and Computation
 \\School of Mathematics and Statistics\\
The University of Western Australia\\
35 Stirling Highway\\
Crawley WA 6009\\
Australia\\
alice.devillers@uwa.edu.au\quad cheryl.praeger@uwa.edu.au
}
\date{today}


\date\today

\maketitle
\begin{abstract}
 We classify the {\it pairwise transitive} 2-designs, that is, 2-designs such that a group of automorphisms is transitive on the following five  sets of ordered pairs: point-pairs, incident point-block pairs, non-incident point-block pairs, intersecting block-pairs and non-intersecting block-pairs.
These 2-designs fall into two classes: the symmetric ones and the quasisymmetric ones. The symmetric examples include the symmetric designs from projective geometry, the 11-point biplane, the Higman-Sims design, and designs of points and quadratic forms on symplectic spaces. The quasisymmetric examples arise from affine geometry and the point-line geometry of projective spaces, as well as several sporadic examples.
\end{abstract}

\section{Introduction}

A design $\cal D$ consists of two sets ${\cal P}$ and ${\cal B}$ of `points' 
and `blocks' respectively, and an incidence relation ${\cal I}\subseteq {\cal P}
\times {\cal B}$; we write  ${\cal D}=({\cal P}, {\cal B}, {\cal I})$.  
The study of designs has a long history, and recurring 
themes are issues of balance and symmetry. Indeed (according to \cite[p. 12]{ACDG}, citing Ahrens~\cite{Ahr}), 
Latin square amulets go back to c.1200, and the study of designs as we have 
formulated them, goes back at least to 1835, when Pl\"ucker, in a study of 
algebraic curves, encountered a Steiner triple system on 9 points (and claimed 
that Steiner systems could only exist if the number of points is congruent to $3\pmod{6}$, a 
conjecture Pl\"ucker later correctly revised to $1$ or $3\pmod{6}$), see \cite[p.12]{ACDG}. 
The mutual importance of groups and designs has also been recognised for decades, for example, 
Witt's discovery of the Steiner systems now known as the Witt designs made their automorphism 
groups much better understood -- these groups are the sporadic simple  Mathieu groups  which had 
been discovered 70 years earlier, see \cite[Chapter IV]{BJL}.

The designs studied in this paper are $2$-designs, that is to say, each 
block is incident with the same number $k$ of points and each pair of 
distinct points has the same number $\lambda$ of incident blocks in 
common. If $v=|{\cal P}|$, we call $\cal D$ a $2$-$(v,k,\lambda)$ design. 
Our aim is to classify all $2$-designs ${\cal D}=({\cal P}, {\cal B}, 
{\cal I})$ which have a strong form of symmetry on pairs 
from ${\cal P}\cup {\cal B}$, defined in the following paragraph. These so-called pairwise transitive designs 
were introduced in \cite{starquo}, where they arose in 
the study  of locally $4$-distance transitive 
graphs. The definition for $2$-designs is a bit simpler than the one given 
in \cite{starquo} since for a $2$-design all point-pairs are incident with 
$\lambda$ common blocks.

An \emph{automorphism} of $\cal D$ is a permutation of ${\cal P}\cup {\cal B}$  
leaving invariant ${\cal P}$ and ${\cal B}$ setwise and preserving incidence.  
A point-block pair is \emph{incident} if it lies in $\cal I$; a block-pair is 
\emph{intersecting} if there is 
a point incident with both blocks. For a subgroup $G$ of the automorphism group $\Aut({\cal D})$ of a $2$-design  
$\cal D$, we say that  $\cal D$ is \emph{$G$-pairwise transitive} if $G$ is transitive 
on the following five (possibly empty) sets of ordered pairs: point-pairs, incident point-block pairs, 
non-incident point-block pairs, intersecting block-pairs and non-intersecting block-pairs. Note that all these pairs are ordered pairs, and that will be the case in the rest of this paper where we will often refer simply to pairs, rather than ordered pairs.
We note that a trivial example is obtained by taking $\cal B$ to be the set of all 
$2$-subsets of a $v$-set $\cal P$ with inclusion as incidence. In the following we 
assume that $\cal D$ is a \emph{non-trivial} $2$-$(v,k,\lambda)$ design in the sense 
that $2< k<v$.  We denote by $\mu$ the number of points incident with the two blocks in an intersecting block-pair.

Transitivity on ordered point-pairs is, in the language of permutation groups,  the property 
of being $2$-transitive on points, and the finite $2$-transitive permutation groups 
are known explicitly as a consequence of Burnside's Theorem and the classification of the finite simple groups 
(see for example, \cite[Sections 7.3 and 7.4]{Cam99}). This classification suggested 
to us the possibility of classifying the pairwise transitive $2$-designs completely, and 
this classification is the aim of our paper.

\begin{theorem}\label{main}
Let $\calD$ be a non-trivial  $2$-$(v,k,\lambda)$ design, and $G\leq \Aut({\cal D})$. Then $\calD$ is $G$-pairwise 
transitive if and only if ${\cal D}$, $G,v,k,\lambda$ are as in one of the rows of Tables \ref{maintable1} and \ref{maintable2}. 
\end{theorem}
   
Note that Table  \ref{maintable2} also gives $\mu$ (we do not list $\mu$ in Table  \ref{maintable1} because $\mu=\lambda $ for symmetric designs). The examples from the two tables are described in detail later in the paper, refer to the column ``Ref.'' in the tables. 
Some of the transitivity conditions in the definition of pairwise transitivity have been studied 
previously for $2$-designs, but not all of them have been imposed at once. For example, for a 
$2$-design $\cal D$ the number of blocks is at least $v=|\cal P|$ and if equality holds then 
$\cal D$ is called a \emph{symmetric $2$-design}. Kantor  \cite{Kan85} applied the classification 
of finite $2$-transitive permutation groups to classify the point 2-transitive symmetric 
$2$-designs.
The non-trivial pairwise transitive symmetric 2-designs coincide with the 2-transitive ones (Lemma \ref{conv-facts}(B)) so most of our effort is focused on the non-symmetric cases. 
Another family of well-studied $2$-$(v,k,\lambda)$ designs is the subfamily with 
$\lambda=1$; these are usually called \emph{linear spaces}. Kantor~\cite{Kan85b} classified the point 
$2$-transitive linear spaces, and a much more general classification was embarked on 5 years 
later: namely the classification of the flag-transitive linear spaces. Flag-transitivity is 
another name for transitivity on incident point-block pairs. This major classification was 
announced in 1990  \cite{BDDKLS} and the last part of the proof was completed in 2003   \cite{S03}; 
the classification leaves open a difficult $1$-dimensional affine case where `complete 
classification is [believed by some to be] hopeless' \cite{Kan93}. On the other hand 
transitivity on non-incident point-block pairs is called antiflag-transitivity. This 
property has been studied to a lesser extent:  Delandtsheer \cite{Del84} classified 
the finite antiflag-transitive linear spaces, and Cameron and Kantor \cite{CamKan79} 
classified the groups of semilinear transformations acting antiflag-transitively on 
some well-known designs. Linear spaces with a group transitive on ordered or unordered pairs of intersecting lines have been studied in \cite{Buek72,Del86,Del92}, while linear spaces with a group transitive on both ordered  pairs of intersecting lines and ordered pairs of disjoint lines are classified by Delandtsheer in \cite{Del84b}.
The assumption of transitivity on ordered pairs of intersecting lines, for a resolvable design, implies transitivity on ordered pairs of distinct parallel classes, that is, 2-transitivity on parallel classes of lines. This property was studied by Czerwinski \cite{Cze}.

In Lemmas \ref{facts} and \ref{sym-qs} we prove that a  non-trivial $G$-pairwise transitive 2-design is such that $G$ acts faithfully with rank 2 or 3 on blocks, in which case the design is symmetric or quasisymmetric respectively. A {\it quasisymmetric design} is a design with exactly two intersection numbers for block-pairs. 
Section \ref{sec:sym} takes care of the former case and is essentially a commentary on Kantor's classification \cite{Kan85}, specifying the groups $G$ and also the examples which are complements of trivial designs.  In the latter case, the action on $\calB$ can either be primitive or imprimitive. The primitive case is treated in Section \ref{sec:prim} and depends on a blend of the classifications of 2-transitive and rank 3 almost simple groups. If $G$ acts imprimitively of rank 3 on $\calB$, then the study divides further according to the type of action on points: affine (Section \ref{sec:affine}) or almost simple (Section \ref{sec:AS}).
Putting all these results together, we obtain Theorem \ref{main}. Note that all computer checks mentioned in this paper use the computer algebra system Magma \cite{magma}.

In \cite{starquo}, we showed that a graph is locally $(G,4)$-distance transitive and has a normal star quotient $K_{1,r}$ for $r\geq 3$ if and only if  its adjacency design is   $G$-pairwise transitive and $N$-nicely affine with $r$ parallel classes of blocks (see Definition \ref{def:nicelyaffine}).
Identifying the $N$-nicely affine designs in Tables \ref{maintable1} and \ref{maintable2} enables us to classify these graphs in the case where the adjacency design is a 2-design.
\begin{corollary}\label{cor}
 Let $\calD$ be a non-trivial  $2$-$(v,k,\lambda)$ design, and $G\leq \Aut({\cal D})$. Then $\calD$ is $G$-pairwise 
transitive and $N$-nicely affine for some non-trivial normal subgroup $N$ of $G$ if and only if ${\cal D}$, $G,v,k,\lambda$ are as
  in one the first three lines of Table \ref{maintable2}.
\end{corollary}

\begin{sidewaystable}
\begin{tabular}{|c|c|p{.8cm}|p{2.3cm}|p{2.3cm}|p{1cm}|c|c|}
\hline
& Design&$v$&$k$&$\lambda$&Ref.& $G$& conditions if any\\
\hline
1&$C_{v,v-1}$&$v$&$v-1$&$v-2$&\ref{complete}&2-transitive subgroup of $\Sym(v)$&\\
\hline
2&$\PG_{d-2}(d-1,q)$ & $\frac{q^d-1}{q-1}$&$\frac{q^{d-1}-1}{q-1}$&$\frac{q^{d-2}-1}{q-1}$ &\ref{PSLdqex} &$\PSL(d,q)\lhd G\leq \PGammaL(d,q)$&$d>2$\\
3& &$15$&$7$&$3$&&$\Alt(7)$&$(d,q)=(4,2)$\\
\hline 
4&$\PG_{d-2}(d-1,q)^c$ & $\frac{q^d-1}{q-1}$&$q^{d-1}$&$q^{d-2}(q-1)$ &comp. of \ref{PSLdqex} &$\PSL(d,q)\lhd G\leq \PGammaL(d,q)$&$d>2$\\
5& &$15$&$8$&$4$&&$\Alt(7)$&$(d,q)=(4,2)$\\
\hline
6&$H(11)$ &$11$&$5$&$2$&\ref{Hadamard11}&$\PSL(2,11)$&\\
\hline
7&$H(11)^c$ &$11$&$6$&$3$&comp. of \ref{Hadamard11}&$\PSL(2,11)$&\\
\hline
8&$D_{176}$ & $176$&$50$&$14$&\ref{Higman}&$\HS$&\\
\hline
9&$D_{176}^c$ & $176$&$126$&$90$&comp. of \ref{Higman}&$\HS$&\\
\hline
10&$S^-(2m)$ & $2^{2m}$&$2^{2m-1}- 2^{m-1}$&$2^{2m-2}- 2^{m-1}$ &\ref{quadr}& $2^{2m}:G_0$ where $\Sp(2m/e,2^e)' \lhd G_0 \leq\Sp(2m,2)$  &$m\geq 2$ and $e\mid m$\\
11&&&&&&$2^{2m}:G_0$ where $\G_2(2^{m/3})'\lhd G_0 \leq\Sp(6,2^{m/3}).m/3$&$m\geq 3$ and $3\mid m$\\
\hline
12&$S^+(2m)$ & $2^{2m}$&$2^{2m-1}+ 2^{m-1}$&$2^{2m-2}+ 2^{m-1}$ &comp. & $2^{2m}:G_0$ where $\Sp(2m/e,2^e)' \lhd G_0 \leq\Sp(2m,2)$  &$m\geq 2$ and $e\mid m$\\
13&&&&&of \ref{quadr}&$2^{2m}:G_0$ where $\G_2(2^{m/3})'\lhd G_0 \leq\Sp(6,2^{m/3}).m/3$&$m\geq 3$ and $3\mid m$\\

\hline
\end{tabular}
\caption{Pairwise transitive symmetric 2-designs}\label{maintable1}

\vspace*{.5cm}
\begin{tabular}{|c|c|p{.8cm}|p{1cm}|p{1cm}|p{1cm}|p{1cm}|c|c|}
\hline
 &Design&$v$&$k$&$\lambda$&$\mu$&Ref.& $G$& conditions if any\\
\hline
1&$\AG(f,q)$&$q^f$&$q^{f-1}$&$\frac{q^{f-1}-1}{q-1}$&$q^{f-2}$&\ref{affine} &$\ASL(f,q)\leq G\leq \AGammaL(f,q)$&\\
2&&$16$&$8$&$7$&$4$&&$2^4:\Alt(7)$&$(f,q)=(4,2)$\\
3&&$16$&$4$&$1$&$1$&&$\AGammaL(1,16)$&$(f,q)=(2,4)$\\
4&&$8$&$4$&$3$&$2$&&$\PSL(2,7)$&$(f,q)=(3,2)$\\
\hline
5&$H(12)$&$12$&$6$&$5$&$3$&\ref{Golay}&$M_{11}$&\\
\hline
6&$\PG_{1}(d-1,q)$&$\frac{q^d-1}{q-1}$&$q+1$&$1$&$1$&\ref{projspace}&$\PSL(d,q)\lhd G\leq \PGammaL(d,q)$&$d\geq 4$\\
\hline
7& $\PG(2,4)^\text{hyp}$& $21$&$6$&$4$&$2$&\ref{hyperovals}&$\PSL(3,4)\lhd G\leq \PSigmaL(3,4)$&\\
\hline
8& ${\cal M}_{22}$&$22$&$6$&$5$&$2$&\ref{Mathieu}&$M_{22}\leq G\leq M_{22}.2$&\\
\hline
\end{tabular}
\caption{Pairwise transitive quasisymmetric 2-designs}\label{maintable2}
\end{sidewaystable}

\section{Notation and preliminary results}

A 2-design, defined above, is actually a particular case of a $t$-design:
a design is called a {\it $t-(v,k,\lambda)$-design} if $|{\cal P}|=v$, 
each block is incident with $k$ points ($t\leq k$), and each $t$-subset of points is incident with exactly $\lambda$ blocks.
A {\it $t$-design} is a  $t-(v,k,\lambda)$-design for some parameters $v,k,\lambda$. 

Let $\calD$ be a $2$-$(v,k,\lambda)$ design. 
We denote by $[p]$ the set of blocks  incident with the point $p$ of $\calD$ and by $[b]$ the set of points incident with the block $b$ of $\calD$. Note that $|[b]|=k$ and we denote $|[p]|$ by $r$ (the number of blocks incident with a point). A block $b$ is called a {\it repeated block} if $[b]=[b']$ for some block $b'\neq b$. 
The \emph{ incidence graph} of  $\calD=({\cal P}, {\cal B}, {\cal I})$ is the graph with vertex-set ${\cal P}\cup {\cal B}$ and $\{x,y\}$ is an edge exactly when  $(x,y)$ or $(y,x)$ is in ${\cal I}$ (less formally: if they form an incident point-block pair). The \emph{point-graph} of  $\calD$ has vertex-set ${\cal P}$ and $\{x,y\}$ is an edge exactly when there is a block of ${\cal B}$ incident with both $x$ and $y$.

 A design $\cal D$ is {\it connected} if its incidence graph  is connected, which is equivalent to its point-graph being connected.
For example, if $\calD$ is a 2-design then its point-graph is a complete graph, and in particular $\calD$ is connected. It also follows that the set of intersecting block pairs of a 2-design is always non-empty.

We say that a $t$-design $\calD$ is {\it non-trivial} if $1\leq t<k<v$, and otherwise $\calD$ is said to be {\it trivial}. For 2-designs this means that  blocks are incident with at least 3 points and not with  all the points.

\begin{lemma}\label{lem:trivial}
Let  $\calD=(\calP,\calB,{\cal I})$ be a  $2$-$(v,k,\lambda)$ design, and $G=\Aut(\calD)$.
Then $\calD$ is trivial and $G$-pairwise transitive if and only if 
either
\begin{description}
 \item[(a)] $k=2$ and $\lambda=1$; or
\item[(b)] $k=v$ and $\lambda=|\calB|$.
\end{description}
\end{lemma}
\proof
First assume that $\calD$ is a trivial and $G$-pairwise transitive 2-design.
Since $\calD$ is trivial, either $k=2$ or $k=v$.
If $k=v$, then all blocks of the design are incident with all the points, so any two points are in all blocks, and  $\lambda=|\calB|$.
If $k=2\neq v$, then $\calP$ contains (at least) three distinct points  $p_1,p_2,p_3$, and since $\calD$ is a 2-design with $k=2$, for each $\{i,j\}\subset \{1,2,3\}$ there is a block $b_{ij}$ such that $[b_{ij}]=\{p_i,p_j\}$. In particular $[b_{12}] \cap [b_{13}]=\{p_1\}$, and since $G$ is transitive on intersecting block-pairs, it follows that no distinct blocks $b,b'$ can satisfy $[b]=[b']$. Thus $\lambda=1$.

By definition of trivial designs, any design satisfying (a) or (b) is trivial. 
In case (b), $G=\Aut(\calD)=\Sym(v)\times \Sym(|\calB|)$, and it is easy to check that  $\calD$ is $G$-pairwise transitive (note that the following sets are empty: non-incident point-block pairs, non-intersecting block-pairs). 
In case (a), $G=\Aut(\calD)=\Sym(v)$, and it is easy to check that  $\calD$, which can simply be seen as a complete graph, is $G$-pairwise transitive.
\qed

We say that a transitive permutation group has {\it rank} $n$ if a point stabiliser has $n$ orbits. Hence a transitive group has rank 2 if and only if it is 2-transitive.

%

\begin{lemma}\label{facts}
Suppose  $\calD$ is a non-trivial $G$-pairwise transitive 2-design.  Then $\calD$ is connected, and $\calD$  and $G$ satisfy the following properties:
\begin{description}
\item[(a)] $G^\calP$ is $2$-transitive;
\item[(b)] Either $G^\calB$ is $2$-transitive, or $G^\calB$ is transitive of rank $3$ and $\calD$ has non-intersecting blocks; 
\item[(c)] $G_p$ has two orbits on $\calB$, for any $p\in\calP$;
\item[(d)] $G_b$ has two orbits on $\calP$, for any $b\in\calB$;
\item[(e)] $\calD$ has no repeated blocks;
\item[(f)] if $p_1,p_2$ are distinct points, then $[p_1]\neq [p_2]$;
\item[(g)] $G$ acts faithfully on  $\calP$;
\item[(h)] $G$ acts faithfully on  $\calB$.
\end{description}
\end{lemma}
\proof
As discussed above, $\calD$ is connected. Part (a)  holds since $G$ is transitive on  point-pairs. 
Since $G$ must be transitive on pairs of intersecting blocks and on pairs of non-intersecting blocks, (b) follows (second part of the statement if there are  non-intersecting blocks, first part otherwise).
Since $\calD$ is non-trivial, $2<k<v$, and so there exist incident point-block pairs, and non-incident point-block pairs. Then, since $\calD$ is $G$-pairwise transitive, $G$ has  exactly  two orbits on $\calP\times\calB$ which implies both (c) and (d).

Suppose $\calD$ has repeated blocks $b_1$, $b_2$, that is $[b_1]=[b_2]$. 
 Since $k\neq v$, there must be a point not incident with $b_1$, say $p$. 
Let $q$ be a point incident with $b_1$. Then there must be a block $b_3$ containing $p$ and $q$ (since $\calD$ is a 2-design). We have that $[b_1]\cap [b_3]\neq \emptyset$, so there must be 
an element of $G$ mapping $(b_1,b_2)$ to $(b_1,b_3)$. Since $[b_1]\neq [b_3]$, this is a contradiction. Hence (e) holds.

Suppose $p_1,p_2$ are distinct points with $[p_1]= [p_2]$. Since $G$ has only one orbit on pairs of points by (a), for any two points $a,b$ we have $[a]= [b]$. By connectedness, all the points are incident
with exactly the same set of blocks, so each block contains all the points, and $\calD$ is trivial, which is a contradiction. Thus (f) holds.

Property (g), respectively (h), follows directly from (e), respectively (f). 
\qed

Thanks to property (e), for non-trivial pairwise transitive designs, we may identify blocks with the sets of points they are incident with, and we will therefore often say that a point is in a block or that a block is a subset of points. Also, because of property (g), we will often  identify $\Aut(\calD)$ with a subgroup of $\Sym(\calP)$.

\begin{lemma}\label{conv-facts}
 Suppose $\calD$ is a non-trivial  2-design, and $G\leq\Aut(\calD)$. 
\begin{description}
 \item[(A)] The following are equivalent:
\begin{enumerate}
 \item[(i)]  $\calD$ is $G$-pairwise transitive,
 \item[(ii)]  Conditions (a), (b), (c) of Lemma \ref{facts} all hold,
\item[(iii)]  Conditions (a), (b), (d) of Lemma \ref{facts} all hold.
\end{enumerate}
\item[(B)] Assume in addition that $\calD$ is symmetric.  Then  $\calD$ is $G$-pairwise transitive if and only if $G^\calP$ is $2$-transitive.
\end{description}
\end{lemma}
\proof
(A)
By Lemma \ref{facts}, (i) implies (ii).
If part (ii) holds, then $G$ has two orbits in its induced action on $\calP\times\calB$ and this in turn implies that condition (d) of Lemma \ref{facts} holds and hence part (iii) holds. 
Finally assume that (iii) holds.
 Condition  (a) of Lemma \ref{facts} implies that $G$ is transitive on point-pairs.  
By (d)  and since $G$ is transitive on points and blocks (by (a) and (b)), $G$ has two orbits on  point-block pairs. Since the sets of  incident point-block pairs and of  non-incident point-block pairs are both non-empty, they must be the two orbits of $G$ on point-block pairs.  
By (b),  $G$ has rank 2 or 3  on $\calB$. In the rank 3 case, $\calD$ has non-intersecting blocks, and the two orbits of $G$ on pairs of distinct blocks must be the pairs of intersecting blocks and the pairs of non-intersecting blocks. If $G$ has rank 2  on $\calB$, then the set of non-intersecting block-pairs is empty and $G$ is transitive on the pairs of intersecting blocks.
It follows that $\calD$ is $G$-pairwise transitive and so (i) holds.

(B)
If $\calD$ is $G$-pairwise transitive then, by part (A), Lemma \ref{facts}(a) holds, which is equivalent to $G^\calP$ being $2$-transitive for a $2$-design $\calD$. 
Conversely suppose that $\calD$ is symmetric and  $G^\calP$ is  $2$-transitive, so Lemma \ref{facts}(a) holds. 
By \cite[Theorem 3.4]{Lander}, $G$ has the same rank considered as a permutation group on points or blocks, so  Lemma \ref{facts}(b) holds. 
Let $p\in\calP$. Since $G$ is 2-transitive on $\calP$, $G_p$ has two orbits on $\calP$. By \cite[Theorem 3.3]{Lander}, $G_p$ also has two orbits on $\calB$, so  Lemma \ref{facts}(c)  holds.
Applying part (A),  $\calD$ is $G$-pairwise transitive.
\qed

For a design $\calD$, we define its {\it complement} ${\cal D}^c$ as the design with the same point set as $\calD$ but with blocks the complements of the blocks of ${\cal D}$ (if we identify each block with the set of points it is incident with). More formally, if  ${\cal D}=({\cal P}, {\cal B}, {\cal I})$, then ${\cal D}^c=({\cal P}, {\cal B}, {\cal I}^c)$, where ${\cal I}^c$ is defined by $(p,b)\in {\cal I}^c$ if and only if  $(p,b)\notin {\cal I}$.
On the other hand, its {\it  dual} is the design ${\cal D}^*=({\cal B}, {\cal P}, {\cal I}^*)$, where  $(b,p)\in {\cal I}^*$ if and only if  $(p,b)\in {\cal I}$.
\begin{lemma}\label{sym-qs}
Let   $\calD$ be a non-trivial  $G$-pairwise transitive $2$-design. 
If $G$ has  rank $2$ on blocks, then $\calD$ is a symmetric design. If $G$ has  rank $3$ on blocks, then $\calD$ is a quasisymmetric design containing disjoint block pairs.
\end{lemma}
\proof
Fischer showed that a non-trivial 2-design has at least as many blocks as points \cite{Fischer}.
Suppose first that $G$ has  rank 2 on blocks. Then any two blocks intersect in a constant number of points, which means that the dual of $\calD$ is also a 2-design. Therefore by Fischer's theorem, $v=b$ and $\calD$ is a symmetric 2-design.

Now assume $G$ has  rank 3 on blocks, so $G$ has two orbits on pairs of distinct blocks, and since $G$ is transitive on intersecting block-pairs and on non-intersecting block-pairs, it follows that there are two intersection sizes for blocks, that is, $\calD$ is a quasisymmetric design. 
\qed

\begin{definition}\label{complete}{\rm
 The {\em complete design} $C_{v,k}$ (where $1<k\leq v$ are integers) is the design with $v$ points such that for each $k$-subset $S$ of points there is exactly one block $b$ with $[b]=S$. In less formal words, each $k$-subset constitutes one block.}
\end{definition}

Since $C_{v,k}$ has no repeated blocks, the action of the automorphism group is faithful on points, and so $\Aut(C_{v,k})\subseteq \Sym(v)$. It follows easily that  $\Aut(C_{v,k})=\Sym(v)$.
By construction, $C_{v,k}$ is a 2-design for any parameters $v$ and $k$.

\begin{lemma}\label{lem-comp}
\begin{description}
 \item[(a)] The complete design $C_{v,k}$ is an $\Sym(v)$-pairwise transitive $2$-design if and only if $k=2$,  $k=v-1$, or $k=v$.
\item[(b)] The complete design $C_{v,k}$ is non-trivial and $G$-pairwise transitive if and only if $k=v-1\geq 3$ and $G$ is a 2-transitive subgroup of $\Aut(C_{v,k})=\Sym(v)$ as in Line 1 of Table \ref{maintable1}. Moreover $C_{v,v-1}$ is a symmetric 2-design.
\end{description}

\end{lemma}
\proof
(a) If $k=v$, then $\calB$ consists of only one block $b$, such that $[b]=\calP$, and $\calD$ is trivially $\Sym(v)$-pairwise transitive.

Next suppose $k=2<v$. We can think of $C_{v,k}$ as the complete graph on $v$ points, with $\calB$ the set of edges. All the required transitivity properties, for $v\geq 4$, follow easily from the fact that $\Sym(v)$ is 4-transitive on the point-set.
If $v=3$ only 3-transitivity is needed (as there are no pairs of non-intersecting blocks).

Now suppose $k=v-1>1$. Each point is incident with all but one block and each block is incident with all but one point.  
If we identify each block with the one point it is not incident with, we see that the automorphism group also acts faithfully as $\Sym(v)$ in its natural action on blocks. 
The design has no pairs of non-intersecting blocks and each point is in only one
non-incident point-block pair. All the required transitivity properties follow easily from the fact that $\Sym(v)$ is 2-transitive on the point-set and on the block-set.

Finally suppose  $C_{v,k}$ is  $\Sym(v)$-pairwise transitive. Then $\Sym(v)$ is transitive on pairs of intersecting blocks, so there can only be one non-zero intersection size for blocks.
If $k<v-1$, then we easily see that there are pairs of blocks intersecting in $k-1$ and $k-2$ points. Thus $k-2$ has to be equal to $0$ (corresponding to non-intersecting blocks).

(b) Let $\calD= C_{v,k}$ be non-trivial and $G$-pairwise transitive. Then $\calD$ is also $\Aut(\calD)$-pairwise transitive, and so by Part (a) $k=2$,  $k=v-1$, or $k=v$. Since $\calD$ is non-trivial, $2<k<v$ and so $k=v-1\geq 3$. Moreover by Lemma \ref{facts}(a), $G^\calP$ is 2-transitive. Since $G$ acts faithfully on points, $G\leq \Sym(v)$. 

Conversely consider $\calD=C_{v,v-1}$, and suppose $G$ is a 2-transitive subgroup of $\Aut(C_{v,k})=\Sym(v)$. Since $\calD$ has $v$ blocks, this design is symmetric, and so by Lemma \ref{conv-facts}(B) it follows that $\calD$ is $G$-pairwise transitive.  Since $k=v-1\geq 3$, $\calD$ is non-trivial. \qed

In a certain sense  $C_{v,v-1}$ is trivial too since it is a $(v-1)-(v,v-1,1)$ design as well as being a 2-design.

\begin{lemma}\label{equiv2trans} Suppose  that $\calD$ is a non-trivial $G$-pairwise transitive $2$-design  and that there exists a $G$-invariant partition $\calS$ of $\calB$ such that $G^\calS$ is $2$-transitive. 
If $G^\calS$ and $G^\calP$ are equivalent, then the  parts of $\calS$ have size $1$, and $\calD$ is the complete design $C_{v,v-1}$.
\end{lemma}
\proof
Let $p\in\calP$. By hypothesis,  $G_p= G_S$ for some part $S\in\calS$. 
Notice that $G_S$ must be transitive on $S$ (since $G$ is transitive on $\calB$), so $S$ is an orbit of $G_S=G_p$ on $\calB$.
By Lemma \ref{facts}(c),  $G_p$ must have two orbits on blocks (blocks incident with $p$, blocks non-incident with $p$). Since $S$ is one $G_S$-orbit in $\calB$, it follows that $G_S$ is transitive on  $\calB\setminus S$. Then, since $G$ preserves incidence, and since $\calD\neq C_{v,v}$, either $[p]=S$ or $[p]= \calB\setminus S$.

If $[p]=S$, then $\calD$ is disconnected, a contradiction. Hence $[p]= \calB\setminus S$, and
  each point  is incident strictly with the blocks not in the corresponding part of $\calB$. Thus, if $b_1,b_2$ are distinct blocks in $S$, then $[b_1]=[b_2]=\calP\setminus\{p\}$,
 contradicting Lemma \ref{facts}(e). 
Therefore the parts of $\calS$ have size 1, and if $S=\{b\}$, then $[b]=\calP\setminus\{p\}$. Hence $\calD=C_{v,v-1}$.
\qed

\section{Pairwise transitive symmetric 2-designs}\label{sec:sym}
In this section, we deal with the case where $\calD$ is non-trivial $G$-pairwise transitive and $G^\calB$ is 2-transitive, so that, by Lemma \ref{sym-qs}, $\calD$ is a symmetric $2$-$(v,k,\lambda)$ design.
By Lemma \ref{conv-facts}(B),  $G^\calP$ is 2-transitive, and the 2-transitive symmetric $2$-$(v,k,\lambda)$ designs $\calD$ with $2<k\leq v-2$ were classified by Kantor \cite{Kan85}. If $k=v-1$ then  $\calD$ is the complete design $C_{v,v-1}$ by Lemma \ref{facts}(e), and is indeed pairwise transitive by Lemma \ref{lem-comp}(b).
We  describe the  designs classified by Kantor and then speicfy explicitly the groups $G$ for which each is $G$-pairwise transitive in Proposition \ref{symmetric2d}.

\begin{example}\label{PSLdqex}{\rm
 Let $\calP$ be the set of points of the projective space $\PG(d-1,q)$ ($d>2$), and let $\calB$ be the set of hyperplanes of  $\PG(d-1,q)$, with inclusion as incidence. 
 Then $\calD=(\calP,\calB,I)$ is a symmetric $2-(\frac{q^d-1}{q-1},\frac{q^{d-1}-1}{q-1},\frac{q^{d-2}-1}{q-1})$ design $\PG_{d-2}(d-1,q)$,  and $\Aut(\calD)=\PGammaL(d,q)$.}
\end{example}

\begin{example}\label{Hadamard11}{\rm
Let $\calD$ be a Hadamard $2-(11,5,2)$ design, shown to be unique up to isomorphism by Todd \cite[Section 3]{Todd33}. 
 Then $\calD$ is a  symmetric $2-(11,5,2)$ design $H(11)$,  and $\Aut(\calD)=\PSL(2,11)$. It is the derived design of a design we will describe in  Example \ref{Golay}.}
\end{example}

The derived design of a given t-design $(\calP,\calB,{\cal I})$ (for simplicity, let us assume it is point-transitive) is the $(t-1)$-design with point-set $\calP \setminus\{p\}$ (where $p$ is a point in $\calP$) and block-set $\{B\setminus\{p\}|B \in \calB, p\in B\}$. 

\begin{example}\label{Higman}{\rm
In the Mathieu-Witt $2-(24,8,5)$ design, fix two points, and consider the two sets of 176 octads that contain precisely one of them. Higman \cite{Hig69} constructed a design, calling the octads in one family ``points" and those in the other family ``quadrics", where a point is incident with a quadric when the octads meet in 0 or 4 symbols. This yields a symmetric $2-(176,50,14)$ design  $D_{176}$, with full automorphism group the Higman--Sims sporadic simple group $\HS$ acting doubly transitively on points and quadrics.
}
\end{example}

\begin{example}\label{quadr}{\rm
Let $V$ be a $2m$-dimensional vector space over $\GF(2)$, with $m\geq 2$. Let $\mathcal{Q}^\epsilon$ ($\epsilon=+$ or $-$) be the set of elliptic (if $\epsilon=+$) or hyperbolic (if $\epsilon=-$) forms on $V$ polarising to a given non-singular alternating form. 
Let $\calP=V$ and $\calB= \mathcal{Q}^+\cup \mathcal{Q}^-$. A vector $v$ is incident with a quadratic form $Q$ of type $\epsilon$ if $(Q(v),\epsilon)=(0,-)$ or $(1,+)$.
The design $\calD=(\calP,\calB,I)$ is denoted by $S^-(2m)$. Its complement is denoted by $S^+(2m)$. 
 Then the design  $S^-(2m)$ is a  symmetric $2-(2^{2m},2^{2m-1}- 2^{m-1},2^{2m-2}- 2^{m-1})$ design,  and $\Aut(\calD)=2^{2m}:\Sp(2m,2)$. }
\end{example}


\begin{proposition}\label{symmetric2d}
 Assume  that $\calD$ is a non-trivial  symmetric $2$-design with $G\leq \Aut({\cal D})$. Then $\calD$ is $G$-pairwise transitive if and only if 
either  ${\cal D}=C_{v,v-1}$ and $G$ is a $2$-transitive subgroup of $\Sym(v)$ (as in Line 1 of Table \ref{maintable1}), or one of  ${\cal D}$ or its complement ${\cal D}^c$ is as in Example \ref{PSLdqex}, \ref{Hadamard11},  \ref{Higman}, or \ref{quadr},  and  $G$ is as in one of the Lines 2, 3, 6, 8, 10, 11 of Table \ref{maintable1}. 
\end{proposition}


\begin{remark}{\rm
 In lines 10, 11 of Table \ref{maintable1}, note that
\begin{enumerate}
\item[(a)]  
$\Sp(2m/e,2^e)' = \Sp(2m/e,2^e)$ unless $(m,e)=(2,1)$, in which case  $\Sp(4,2)' \cong \Alt(6)$ and $\Sp(4,2)\cong \Sym(6)$;
\item[(b)] 
$\G_2(2^{m/3})'=\G_2(2^{m/3})$ unless $m=3$, in which case  $\G_2(2)' \cong \PSU(3,3)$ and $\G_2(2) \cong \PSU(3,3).2$.
\end{enumerate}}
\end{remark}

\proof
Assume $\calD$ is $G$-pairwise transitive and $v=k-1$. Then, as mentioned above,  $G^\calP$ is 2-transitive, and  ${\cal D}=C_{v,v-1}$, so  $G$ is a $2$-transitive subgroup of $\Sym(v)$. 
Conversely, if ${\cal D}=C_{v,v-1}$ and  $G$ is a $2$-transitive subgroup of $\Sym(v)$, then $\calD$ is $G$-pairwise transitive  by Lemma \ref{lem-comp}(b)
Thus we may assume, since $\calD$ is non-trivial, that $2<k\leq v-2$.
Using the well-known equality  \cite[p.3]{Lander} $(v-1)\lambda=k(k-1)$, it is easy to show that $k\neq v/2$ and $k\neq v-2$.
If $k>v/2$ then by Kantor's discussion in the first paragraph of \cite{Kan69}, $\calD^c$ is a symmetric 2-design with block-size $v-k$ strictly less than $v/2$.
So, up to taking the complement, we can assume that $2<k<v/2$.
Thus, by Kantor \cite{Kan85}, $\calD$  is  one of the four examples given in the statement.  
These examples are described in detail in Lander \cite[p.84-89]{Lander}, where the full automorphism groups are described (except for Example \ref{PSLdqex}, but the full group of that design is well known to be $\PGammaL(d,q)$). See also \cite{Kan75} for s description of Example \ref{quadr}.
To find the groups $G$ for which $\calD$ is $G$-pairwise transitive, we only have to find the 2-transitive subgroups of $\Aut(\calD)$, in each case, which yields the groups on the Lines 2, 3, 6, 8, 10, 11 of Table \ref{maintable1}.  For the case of Lines 10-11, $G=2^{2m}:G_0$ with $G_0$ a transitive subgroup of $\Sp(2m,2)$, we refer to Appendix 1 in \cite{Lie87} or to \cite[p.68]{Kan85b}.
%
\qed

 \section{Quasisymmetric pairwise transitive 2-designs}
In this long final section, we deal with the case where $\calD$ is pairwise transitive and $G^\calB$ has rank 3, so that by Lemma \ref{sym-qs}, $\calD$ is a quasisymmetric design containing disjoint block pairs. The analysis splits naturally into the cases where $G^\calB$ is primitive and imprimitive, and the latter case subdivides again into the cases where $G^\calP$  is affine or almost simple 2-transitive (see Lemma \ref{facts}). We treat these subcases in separate subsections below, but first we prove a general lemma.

Recall that $r$ denotes the number of blocks incident with a point.

\begin{lemma}\label{lem:orbitsizes}
 Suppose that $\calD$ is a non-trivial $G$-pairwise transitive $2$-design, and $G^\calB$ has rank $3$. Then $3\leq k\leq v/2$, $3\leq r\leq |\calB|/2$, and for $b\in\calB$, $G_b$ has an orbit on $\calP$ of size $k$.
\end{lemma}
\proof
Since $\calD$ is non-trivial, $k\geq 3$, and since  $G^\calB$ has rank 3, the set of non-intersecting blocks is non-empty, thus $k\leq v/2$.  
Since $|\calB|\cdot k=v\cdot r$, we have $r\leq |\calB|/2$. For each point $p$, there are at most $r(k-1)$ other points lying in blocks incident with $p$. Since $\calD$ is a $2$-design, each point distinct from $p$ lies in one of these blocks and hence $v-1 \leq r(k-1) \leq r(v-2)/2$, so $r\geq 2(v-1)/(v-2)>2$.
By Lemma \ref{facts}(d), for a block $b\in\calB$, $G_b$ has two orbits on $\calP$, namely the set $[b]$ of points incident with $b$, and the rest $\calP\setminus [b]$. The last statement follows.
\qed



\subsection{ $G$ is primitive of rank 3 on blocks }\label{sec:prim}
First we observe that the group $G$ must be almost simple and we give three examples. Our strategy is then to compare the possibilities for the 2-transitive $G^\calP$ and rank 3 $G^\calB$, each isomorphic to $G$, using the classifications of such groups.
\begin{lemma}\label{primblocksAS}
 Suppose  $\calD$  is a non-trivial $G$-pairwise transitive $2$-design, $G^\calB$ has rank $3$ and is primitive. 
Then $G$ is almost simple.
\end{lemma}

\proof
Since $G^\calP$ has rank 2, $G$ is of  either affine or of almost simple type by Burnside's Theorem \cite[p.101]{Cam99}. 
Suppose $G$ is affine.
Let $T=\soc(G)$ be its translation subgroup, which can be identified with $\calP$. 
Since $T\lhd G$ and $G^\calB\cong G$ is primitive, the group $T^\calB\cong T$ is transitive, and hence $|\calB|$ divides $|T|=v$. 
By Fisher's inequality \cite{Fischer}, it follows that $v=|\calB|$, that is $\calD$ is symmetric, contradicting Lemma \ref{sym-qs}.
 \qed

\begin{example}\label{projspace}{\rm
Let $\calP=\PG(d-1,q)$ for $d\geq 4$, $\calB$ be the set of lines of this projective space, with inclusion as incidence.
Then $\calD=(\calP,\calB,I)$ forms a 
$2-(\frac{q^d-1}{q-1},q+1,1)$ design  $\PG_{1}(d-1,q)$. Blocks intersect in 0 or 1 point,  and $\Aut(\calD)=P\Gamma L(d,q)$.
}
 \end{example}

\begin{example}\label{hyperovals}{\rm
Let $\calP=\PG(2,4)$, and let $\calB$ be one of the three orbits of $56$ hyperovals under $\PSL(3,4)$, with inclusion as incidence. 
 Then $\calD=(\calP,\calB,I)$ forms a 
$2-(21,6,4)$ design  $\PG(2,4)^\text{hyp}$. Blocks intersect in 0 or 2 points, and $\Aut(\calD)=P\Sigma L(3,4)$.
}
 \end{example}

\begin{example}\label{Mathieu}{\rm
Let $\calD$ be the Mathieu-Witt design ${\cal M}_{22}$. This is a $3-(22,6,1)$ design and its blocks are often called hexads.  
In particular $\calD$ is a $2-(22,6,5)$ design. 
Blocks intersect in 0 or 2 points,  and $\Aut(\calD)=M_{22}.2$.
} \end{example}

Note that Example \ref{hyperovals} is the residual design of Example \ref{Mathieu}. A {\it residual design} of a  $t$-design $(\calP,\calB,{\cal I})$ (for simplicity, let us assume it is point-transitive) is the $(t-1)$-design with point-set $\calP \setminus\{p\}$ (where $p$ is a point in $\calP$) and block-set $\{B|B \in \calB, p\notin B\}$. 

The number of points incident with two intersecting blocks, the parameters of the design, and the full automorphism group for the three previous examples are 
 well-known.

\begin{theorem}\label{blockaction-rk3-prim}
Assume  $\calD$ is a non-trivial  $2$-design and assume $G\leq \Aut({\cal D})$ is such that  $G^\calB$ has rank $3$ and  is primitive.
Then $\calD$ is $G$-pairwise transitive if and only if  ${\cal D}$ is as in Example \ref{projspace}, \ref{hyperovals},  or \ref{Mathieu},  and  $G$ is as in one of the Lines  6, 7, 8 of Table \ref{maintable2}, respectively. 
\end{theorem}

As mentioned above, $G\cong G^\calP\cong G^\calB$ must arise in the classifications of almost simple groups which are 2-transitive and which are rank 3.
 We need to consider all  pairs of such actions. 
However if up to permutational isomorphism there is only one 2-transitive action and the rank 3 actions are  fused by outer automorphisms (or vice versa), then it is enough to consider just one of the pairs of actions. The list of almost simple 2-transitive groups can be found in \cite{Cam99} and the list of almost simple rank 3 groups follows from Bannai \cite{Bannai}, Kantor and Liebler \cite{KanLie}, and Liebeck and Saxl \cite{LieSaxl}. Comparing the two lists, we list the relevant groups in Tables \ref{prim2-3} and \ref{prim2-3b} (at the end of the paper).
 Column 1 gives $\soc(G)$ (or $G$ in a few cases), Column 2/3 gives the degree of the 2-transitive/rank 3 action (if possible with a description), Column 4 gives the orbit sizes of $G_b$ on $\calP$ (these numbers can often be deduced from our knowledge of the actions, or otherwise using a computer), the entry in Column 5 relates to the proof. 

\medskip
\proof 
Assume  $\calD$ is $G$-pairwise transitive. From the above discussion, one of the lines of Table \ref{prim2-3} or \ref{prim2-3b} holds for $G$. 
The entry in
 Column 5 indicates how we prove that either the group does not provide an example, or corresponds to one of the above examples. The following comments indicate how various cases are dealt with.

\begin{sidewaystable}[H]
\centering

\begin{tabular}{|p{3.5cm}|p{6cm}|p{5.8cm}|p{3.3cm}|p{3.5cm}|}
\hline
$\soc(G)$ (or $G$ if specified) & degree of 2-transitive action &degree of rank 3 action&orbit sizes of $G_b$ on $\calP$ &proof \\
\hline
$\Alt(n)$, $n\geq 5$&$n$ &$\binom{n}{2}$ (on pairs)&$2,n-2$&Lemma \ref{lem:orbitsizes}\\

$\Alt(9)$ &$9$ &$120$ (on non-singular points for $Q^+(8,2)$)&$9 $& Lemma \ref{facts}(d)\\
$\Alt({10})$ &$10$ &$126$ (on $5|5$-partitions of a $10$-set)&$10 $& Lemma \ref{facts}(d)\\
$\Alt(6)$ &$6$ (on totals) &$15$ (on pairs)&$6 $& Lemma \ref{facts}(d)\\
$\Alt({8})$ &$8$ &$35$ (on lines of $\PG(3,2)$/ $4|4$-partitions of a $8$-set)&$8 $& Lemma \ref{facts}(d)\\

$\PSL(d,q)$, $d\geq 4$ &$\frac{q^d-1}{q-1}$ (on points)&$\frac{(q^d-1)(q^{d-1}-1)}{(q-1)^2(q+1)^2}$(on lines)&$q+1,\frac{q^d-q^2}{q-1}$& Example \ref{projspace}\\
$\PSL(d,q)$, $d\geq 5$ &$\frac{q^d-1}{q-1}$ (on points)&$\frac{(q^d-1)(q^{d-1}-1)}{(q-1)^2(q+1)^2}$(on codimension 2 subspaces)&$\frac{q^{d-2}-1}{q-1},\frac{q^{d}-q^{d-2}}{q-1}$& no disjoint blocks\\

$\PSL(2,5)\cong \Alt(5)$ &$6$ (on points)&$10$ (on pairs of 5-set)&$3,3$&no disjoint blocks\\
$\PSL(2,9)\cong \Alt(6)$ &$10$ (on points)&$15$ (on pairs of 6-set)&$4,6$&no disjoint blocks\\
$\PSL(4,2)\cong \Alt(8)$ &$15$ (on points)&$28$ (on pairs of 8-set)&$15$& Lemma \ref{facts}(d)\\

$G=\PGammaL(2,8)$ &$9$ (on points)&$36$ (on pairs of points)&$2,7$&Lemma \ref{lem:orbitsizes}\\

$\PSL(3,4)$ &$21$ (on points)&$56$ (on an orbit of hyperovals)&$6,15$& Example \ref{hyperovals}\\

$\Sp(2d,2)$, $d\geq 3$&$2^{2d-1}+2^{d-1}$ (on hyperbolic forms)&$2^{2d}-1$ (on non-zero vectors)& $2^{2d-2},2^{2d-2}+2^{d-1}$ &no disjoint blocks\\
$\Sp(2d,2)$, $d\geq 3$&$2^{2d-1}-2^{d-1}$ (on elliptic forms)&$2^{2d}-1$ (on non-zero vectors)&$2^{2d-2}-2^{d-1} ,2^{2d-2}$&no disjoint blocks\\

$\PSU(3,5)$&$126$ (on singular points)& $50$&$126 $ & Lemma \ref{facts}(d)\\
$G=\PGammaU(3,3)$&$28$ (on singular points)&$36$&$28$& Lemma \ref{facts}(d)\\
\hline
\end{tabular}
\caption{Alternating and classical groups with both a 2-transitive action and a rank 3 action}\label{prim2-3}
\end{sidewaystable}

\begin{sidewaystable}[H]
\centering

\begin{tabular}{|p{2.7cm}|p{5cm}|p{5.5cm}|p{2cm}|p{3.5cm}|}
\hline
$\soc(G)$ (or $G$ if specified) & degree of 2-transitive action &degree of rank 3 action&orbit sizes of $G_b$ on $\calP$ &proof \\
\hline
$M_{11}$&$11$&$55$ (pairs of 11-set)&$2,9$&Lemma \ref{lem:orbitsizes}\\
$M_{11}$&$12$&$55$ (pairs of 11-set)&$12$& Lemma \ref{facts}(d)\\

$M_{12}$&$12$ (on points)&$66$ (pairs of 12-set)&$2,10$&Lemma \ref{lem:orbitsizes}\\
$M_{12}$&$12$ (on totals)&$66$ (pairs of 12-set)&$12$& Lemma \ref{facts}(d)\\

$\Alt(7)$&$15$ (on points of $\PG(3,2)$)&$21$ (on pairs of 7-set)&$15$& Lemma \ref{facts}(d)\\

$M_{22}$&$22$&$77$ (on hexads)&$6,16$& Example \ref{Mathieu}\\
$M_{22}$&$22$&$176$ (on heptads)&$7,15$&no disjoint blocks\\

$M_{23}$&$23$&$253$ (on pairs)&$2,21$&Corollary \ref{lem:orbitsizes}\\
$M_{23}$&$23$&$253$ (on heptads)&$7,16$&no disjoint blocks\\

$M_{24}$&$24$&$276$ (on pairs)&$2,22$&Lemma \ref{lem:orbitsizes}\\
$M_{24}$&$24$&$1288$ (on pairs of dodecads)&$24$& Lemma \ref{facts}(d)\\

$G=\PGammaL(2,8)$  &$28$&$36$ (on pairs of 9-set)&$7,21$&no disjoint blocks\\

$HS$&$176$&$100$&$176$& Lemma \ref{facts}(d)\\
\hline
\end{tabular}
\caption{Sporadic almost simple groups with both a 2-transitive action and a rank 3 action}\label{prim2-3b}
\end{sidewaystable}
\begin{itemize}
 \item The groups with ``Lemma \ref{facts}(d)'' in the last column cannot be the automorphism group of a design satisfying the required properties, as $G_b$ has only one orbit on $\calP$, contradicting Lemma \ref{facts}(d).

\item The groups with ``Lemma \ref{lem:orbitsizes}'' in the last column cannot be the automorphism group of a design  satisfying the required properties, as $G_b$ has no orbit of size between $3$ and $v/2$ on $\calP$, contradicting Lemma \ref{lem:orbitsizes}.
In all other cases, by Lemma \ref{lem:orbitsizes}, a given block $b$ must be incident to the points in $\calP$ in the smaller $G_b$-orbit on $\calP$. 
In each case, we call $\calD$ this design.

\item The groups with ``no disjoint blocks'' in the last column cannot be the automorphism group of a design  satisfying the required properties, as we will show there are no disjoint blocks, and instead there are two orbits on pairs of intersecting blocks. 
\end{itemize}
Consider the case of $G=\PSL(d,q)$ in its rank 3 action on codimension 2 subspaces, where $d\geq 5$ (if $d=4$, this action is the same as the action on lines, as in the previous line of the table). 
Since the number of points contained in a codimension 2 subspace is smaller than the number of points outside it, $\calD$ is the design of projective points/codimension 2 subspaces, with inclusion as incidence. Two blocks intersect in a subspace of codimension 3 or 4, each of which contains a 1-dimensional subspace (projective point). 

Consider the case of $G=\Sp(2d,2)$ ($d\geq 3$). Then $\calP$ is the set of quadratic forms of type $\epsilon$ (where $\epsilon=\pm $) polarising to the symplectic form $(\,,\,)$ preserved by $G$, and $\calB$ is the set of non-zero vectors of $V(2d,2)$. Let $p$ be a point (so a quadratic form). Then $p$ has $2^{2d-1}+\epsilon 2^{d-1} -1$ singular non-zero vectors and  $2^{2d-1}-\epsilon 2^{d-1}$ non-singular vectors, and these form the two orbits of $G_p$ on $\calB$. By  Lemma \ref{lem:orbitsizes}, the smallest one of these orbits will give us  incidence. For $\epsilon=+$ (hyperbolic forms), the set of non-singular vectors is smaller, so a point $p$ and a block $b$ are incident in $\calD$ if and only if $b$ is non-singular for the hyperbolic form $p$.	
For $\epsilon=-$ (elliptic forms), the set of non-zero singular vectors is smaller, so a point $p$ and a block $b$ are incident in $\calD$ if and only if $b$ is singular for the elliptic form $p$.

We claim that, in both cases, there are no disjoint blocks, that is for any two blocks $x,y$, there exists a point incident with both. In other words for any two non-zero vectors $x,y$ of $V(2d,q)$, there exists a hyperbolic (resp. elliptic) form polarising to $(,)$ such that $x,y$ are both non-singular (resp. singular) with respect to this form. As $G$ has two orbits on pairs of blocks, it is sufficient to consider two pairs $x,y$: one where $(x,y)=0$ and one where $(x,y)=1$. 
If $x$ is the vector $(x_1,x_2,\ldots,x_{2d})$ (and similarly for $y$), we can pick without loss of generality the bilinear form to be $(x,y)=\sum_{i=1}^d (x_{2i-1}y_{2i}+x_{2i}y_{2i-1})$.

Let $x=(1,1,0,\ldots,0)$, $y=(0,0,1,1,0,\ldots,0)$ and $y'=(1,0,1,1,0,\ldots,0)$ be blocks. Notice that  $(x,y)=0$ and  $(x,y')=1$.
Then $\phi$ defined by $\phi(x)=\sum_{i=1}^d x_{2i-1}x_{2i}$ is a hyperbolic form polarising to $(,)$ for which  $x,y, y'$ are all non-singular.

Let $x=(0,0,1,0,\ldots,0)$, $y=(0,0,0,0,1,0,\ldots,0)$ and $y'=(0,0,0,1,0,\ldots,0)$ be blocks. Notice that  $(x,y)=0$ and  $(x,y')=1$.
Then $\phi$ defined by $\phi(x)=x_1^2+x_2^2+\sum_{i=1}^d x_{2i-1}x_{2i}$ is an elliptic form polarising to $(,)$ \cite[Lemma 2.5.2(ii) and Prop. 2.5.3 (ii)]{KleiLie} for which  $x,y, y'$ are all singular.

Thus the claim is proved.

All the other lines with ``no disjoint blocks'' in the last column  are easily dealt with by computer  (that is, using the computer system Magma), leaving us with precisely three cases, namely $\PG_1(d-1,q)$ with $\soc(G)=\PSL(d,q)$, $\PG(2,4)^\text{hyp}$ with  $\soc(G)=\PSL(3,4)$,  and ${\cal M}_{22}$ with $\soc(G)=M_{22}$. If we add the knowledge of the full automorphism group, then ${\cal D}$ and $G$ are as in the statement.

Conversely, assume that  ${\cal D}$ and $G$ are as in Example \ref{projspace}, \ref{hyperovals},  or \ref{Mathieu},  and  $G$ is as in Lines  6, 7, 8 of Table \ref{maintable2}, respectively. 
In each case, $G^\calP$ has rank 2 (so condition (a) of Lemma \ref{facts} holds) and $G^\calB$ has rank 3 and is primitive.  
We check for each case  that a block stabiliser in $G$ has two orbits on  $\calP$ (see Tables \ref{prim2-3} and \ref{prim2-3b} for sizes):
\begin{description}
 \item[$\PG_1(d-1,q)$:] the stabiliser of a line is transitive on the points on the line, and on the points outside the line in a projective space.
\item[$\PG(2,4)^\text{hyp}$:] the stabiliser of a hyperoval is transitive on the 6 points of the hyperoval, and on the 15 points outside the hyperoval (this is well known and can be deduced from the character table of $\PSL(3,4)$, see the Atlas \cite[p. 23]{atlas}).
\item[${\cal M}_{22}$:] the stabiliser of a hexad is transitive on the 6 points of the hexad, and on the 16 points outside the hexad (this is well known and can be deduced from the character table of $M_{22}$, see the Atlas \cite[p. 39]{atlas}).
\end{description}
 Hence  condition (d) of Lemma \ref{facts} holds. 
To verify condition (b) of Lemma \ref{facts}, we need to check that there are non-intersecting blocks. This is obvious for $\PG_1(d-1,q)$, as $d\geq 4$ so there are disjoint lines. For $\PG(2,4)^\text{hyp}$ and  ${\cal M}_{22}$, we checked this by computer. Thus condition (b) of Lemma \ref{facts} holds.
By Lemma \ref{conv-facts}, $\calD$ is $G$-pairwise transitive. 
\qed

\subsection{$G$  is imprimitive of rank 3 on blocks}
In this case, the group $G$ is affine or almost simple. 
Recall from \cite{starquo} the following definition (which is valid for any $t$-design not just $2$-designs). 
\begin{definition}\label{def:nicelyaffine}
{\rm Let ${\cal D}=({\cal P}, {\cal B}, {\cal I})$ be a design and $N$ be an automorphism group of ${\cal D}$. Then  ${\cal D}$ is called {\it $N$-nicely affine} if $N$ is transitive on ${\cal P}$, and there is a constant $\mu$ such that distinct blocks  are incident with exactly $\mu$ common points if they are in different $N$-orbits and are disjoint if they are in the same $N$-orbit. 
}
\end{definition}

It follows that the $N$-orbits on the blocks of an $N$-nicely affine are parallel classes.

\begin{lemma}\label{imprimcase}
Suppose  $\calD$  is a non-trivial $G$-pairwise transitive 2-design, and $G^\calB$ has rank $3$.
If $G$ is imprimitive on blocks, then one of the following holds:
\begin{description}
 \item[(1)]  $G$ is $2$-transitive of affine type on $\calP$, and $\calD$ is $N$-nicely affine for $N$ the translation subgroup of $G$;
\item[(2)] $G$  is $2$-transitive of almost simple type on $\calP$ and is quasiprimitive on $\calB$.
\end{description}

\end{lemma}
\proof 
The group $G$ admits a unique nontrivial system of imprimitivity $\calS$ on $\calB$ (see, for example \cite{rk3}).
The equivalence relation given by $\calS$ corresponds to a pair of blocks  being either intersecting, or  non-intersecting.
In the former case, blocks in different parts of $\calS$ are disjoint, and it follows that the design is not connected, a  contradiction. 
Hence distinct blocks in the same part of $\calS$ are disjoint.
Let $M$ be the kernel of the action of $G$ on the system of imprimitivity. 

Suppose $M\neq 1$. 
Then since $G\cong G^\calP$ is 2-transitive,  it follows that $M$ is transitive on points, and so the blocks in each $M$-orbit in $\calB$ 
cover the point-set and form a block of imprimitivity for  $G^\calB$. The uniqueness of $\calS$ implies that each of these blocks is exactly a part of $\calS$. In other words, 
each part of $\calS$ is a parallel class of blocks.  
Since $G$ is transitive on pairs of non-intersecting blocks, it follows that pairs of blocks from different $M$-orbits are intersecting, and since $G$ is transitive on intersecting block-pairs, we conclude that $\calD$ is $M$-nicely affine.

 Since $G\cong G^\calP$ is 2-transitive, the group $G$ is either almost simple or affine.
Let $N=\soc(G)$.  
Then $N$ is the unique minimal normal subgroup of $G$ and $N^\calP$ is transitive. By uniqueness $N\leq M$ and in particular $N$ is intransitive on $\calB$.
Since $G$ admits a unique nontrivial system of imprimitivity on $\calB$, the set of $N$-orbits on $\calB$ must be $\calS$, and so $\calD$ is $N$-nicely affine.
In particular, $N$ preserves  the partition of $\calP$ given by any part of  $\calS$, and so $N$ is imprimitive on points.

By Wielandt \cite[Exercise 12.4]{Wie}, drawing together  \cite[Exercise 12.3, Theorems 12.1, 10.4, 5.1, 11.3]{Wie}, every non-regular minimal normal subgroup of a 2-transitive group is primitive. Hence  $N$ is regular on $\calP$. 
%
By Wielandt \cite[Theorem 11.3(a)]{Wie}, since $G^\calP$ is 2-transitive and $N^\calP$ is a regular normal subgroup, it follows that $N$ is elementary abelian
and hence $G\cong G^\calP$ is of affine type, so (1) holds.

Now suppose $M=1$, that is,  $G$  acts faithfully on the imprimitivity system on blocks. By \cite[Corollaries 2.4 and 2.5]{rk3} $G$ is almost simple and quasiprimitive, so (2) holds.\qed

\subsubsection{$G$ is affine on points}\label{sec:affine}
We study in this section case (1) of Lemma \ref{imprimcase}, that is, $N\lhd G\leq \AGL(d,p)$, and $G$ is 2-transitive of affine type on $\calP$ with translation subgroup $N$ (where $p$ is a prime).
The affine geometry gives a natural example.

\begin{example}\label{affine}
{\rm
 Let $\calP$ be the set of points of the affine space $\AG(f,q)$ ($f>1$), and let $\calB$ be the set of affine hyperplanes of  $\AG(f,q)$, with incidence given by inclusion.
 Then $\calD=(\calP,\calB,I)$ is a $2-(q^f,q^{f-1},\frac{q^{f-1}-1}{q-1})$ design, denoted with a little abuse of notation by $\AG(f,q)$,  and $\Aut(\calD)=\AGammaL(f,q)$. Blocks intersect in $0$ or $q^{f-2}$ points.}
 \end{example}

\begin{theorem}\label{rank3impaffine}
Let  $\calD$ be a non-trivial   $2$-design and let $G\leq \Aut({\cal D})$ be such that  $G^\calB$ has rank $3$ and  is imprimitive, and $G$  is $2$-transitive of affine type on $\calP$.
Then $\calD$ is $G$-pairwise transitive if and only if 
 $\calD$ is the design $AG(f,q)$ described in Example \ref{affine} and  
 $G$ is as in one of Lines  1, 2, 3 of Table \ref{maintable2}.

\end{theorem}

In \cite[Theorem 1.8]{starquo}, we proved that an $N$-nicely affine, $G$-pairwise transitive design, such that $G^\calP$ is 2-transitive and $N$ has at least 3 block orbits, can be obtained from the following construction. The construction has been modified from \cite{starquo} to take into account the condition that $\calD$ is a 2-design.

\begin{construction}\label{constr-regN}
Let $V=V(d,p)$ be a vector space with  translation group $N$, and let  $G=N.G_0\leq \AGL(d,p)$, where $p$ is a prime and $G_0\leq \GL(d,p)$. 
Assume that all the following conditions hold.
 \begin{description}
 \item[(a)] $G_0$ is transitive on $V\setminus\{0\}$;
\item[(b)] there exists a $G_0$-orbit $\calM=\{M_1,\ldots,M_r\}$ ($r\geq 3$) of subspaces of $V$ such that $G_0^{\calM}$ is $2$-transitive (not necessarily faithful);
\item[(c)] $V=M_1+M_2$;
\item[(d)] the stabiliser $(G_0)_{M_1}$ acts transitively on the nontrivial elements of  $V/M_1$;
\end{description} 
Define the design ${\cal D}=(V, \cup_{i=1}^r V/M_i, {\cal I})$ with  incidence $\cal I$ given by inclusion.
\end{construction}

Note that $V/M_i$ denotes the sets of cosets of $M_i$ in $V$.
In our analysis we use the following information about $\Gamma L(1,p^d)$.

\begin{lemma}\cite[Lemma 4.1]{Foulser} and \cite[Lemma 4.7]{LLP} \label{GaLcriterion}
 Let $\FF=\GF(p^d)$ and $\epsilon$ be a primitive element of $\FF$. Then  $\Gamma L(1,p^d)=\langle \tau,\sigma \rangle$, where $\tau:\FF\rightarrow\FF: x\mapsto \epsilon x$ and  $\sigma:\FF\rightarrow\FF: x\mapsto  x^p$. Let $G_0\leq \Gamma L(1,p^d)$. Then there exist  unique integers $i,j,t$ such that $G_0=\langle \tau^i,\tau^j\sigma^t  \rangle$ and the following all hold:
\begin{itemize}
 \item[(1)] $i>0$ and $i\div p^d-1$;
\item[(2)] $t>0$ and $t\div d$;
\item[(3)] $0\leq j<i$ and $i\div j\frac{p^d-1}{p^t-1}$.
\end{itemize}
Moreover $G_0$ is transitive on $V(d,p)\setminus\{0\}$ if and only if either  $i=1$ or the following two conditions are both satisfied:

(A) $j>0$ and $i\div j\frac{p^{it}-1}{p^t-1}$; \hspace{2cm}
(B)  if $1< k<i$ then $i\notdiv j\frac{p^{kt}-1}{p^t-1}$.

\end{lemma}
A generating set for $G_0$ as in Lemma \ref{GaLcriterion} is said to be in {\it standard form}.

\medskip
\noindent{\it Proof of Theorem} \ref{rank3impaffine}.

Assume  $\calD$ is a non-trivial $G$-pairwise transitive  2-design and $G\leq \Aut({\cal D})$ is such that  $G^\calB$ has rank 3 and  is imprimitive, and $G$ is 2-transitive of affine type on $\calP$.
Then  $N\lhd G\leq \AGL(d,p)$  with translation subgroup $N$, where $p$ is a prime, and  $N$ acts regularly on $\calP$.

By  Lemma \ref{imprimcase}, $\calD$ is $N$-nicely affine.  The number of parallel classes of blocks is equal to the number of blocks incident with each point, denoted by $r$. By Lemma \ref{lem:orbitsizes}, we have $r\geq 3$.
It follows that $\calD$ satisfies all the assumptions of \cite[Theorem 1.8]{starquo}, and so $\calD$ can be obtained via  Construction \ref{constr-regN}. 
Thus $G_0$ has a subgroup $H$ of index $r$ such that the  action of $G_0$ on the right cosets of $H$ is 2-transitive. The subgroup $H$ is the stabiliser in $G_0$ of a subspace $M_1$ of $\mathbb{F}_p$-dimension at least $d/2$ (by condition (c)).

The finite 2-transitive groups of affine type have been classified (see for instance \cite{Cam99}):
they consist of 3 infinite families and 7 sporadic cases.
It can easily be checked by computer that out of the 7 sporadic cases, the only case where there is a subgroup $H$ as described is:
$V=V(4,2)$, $G=2^4:\Alt(7)$, $r=15$ and  $\calM$ is the set of all hyperplanes of $V$. Then the blocks are all the cosets of the hyperplanes of $V$, that is, all the affine hyperplanes af $\AG(4,2)$ and Line   2 of Table \ref{maintable2} holds. 
We now look at the 3 infinite families.

\medskip

1. $\SL(f,q)\leq G_0\leq \Gamma L(f,q)$, where $f=d/e$ and $q=p^e$ for some divisor $e$ of $d$ (notice that $p^d=q^f$).

Assume first that $f\geq 2$.  Then  $\SL(f,q)$ is not contained in the kernel of the action of $G_0$ on $\calM$, since it is transitive on the non-zero vectors and so cannot leave a proper subspace invariant. 
Thus either  the 2-transitive group induced by $G_0$ on the $H$-cosets is almost simple with socle $\PSL(f,q)$, or  $(f,q)=(2,2)$ or $(2,3)$. 
Consider first the 2-transitive actions of $G_0$ of degree 
$r=\frac{p^d-1}{q-1}$ (there are  two such actions if $f\geq 3$). The stabiliser of an element in one of these actions stabilises either a $1$-space or a hyperplane of $V(f,q)$. 
Since the subspaces in $\calM$ have $\mathbb{F}_p$-dimension at least $d/2$, the action must be on hyperplanes. Hence in this case $M_1$ is a hyperplane and its orbit under $G_0$ consists of all hyperplanes  of  $V(f,q)$, and Line 1 of Table \ref{maintable2} holds. It is easily checked that the conditions of Construction \ref{constr-regN} are satisfied.
There are a few other possibilities for $2$-transitive actions of $G_0$ if $(f,q)\neq(2,2)$ or $(2,3)$, and they are summarised in the table below (most of them coming from exceptional isomorphisms of $\PSL(f,q)$). It is easily checked by computer that  the stabiliser of an element in each of these actions does not stabilise a subspace.
Finally, assume $(f,q)=(2,2)$ or $(2,3)$. Then, apart from the natural actions of degree $q+1$ already treated, the only possibility is $G_0=\GL(2,3)$ which has a $2$-transitive action of degree $3$. However a stabiliser in that action does not stabilise any subspace.

\begin{center}
\begin{tabular}{|c|c|c||c|c|c|}
 \hline 
$(f,q)$&$\soc(G_0/Z)\cong$&  Degree& $(f,q)$& $\soc(G_0/Z)$ &  Degree\\
\hline
$(2,5)$&$\Alt(5)$& $5$&$(4,2)$&$\Alt(8)$&$8$\\
$(2,7)$&$\PSL(3,2)$&$7$&$(2,8)$&$R(3)'$&$28$\\
$(3,2)$&$\PSL(2,7)$ &$8$ &$(2,11)$&$\PSL(2,11)$& $11$ \\
$(2,9)$&$\Alt(6)$&$6$ &&&\\
\hline
\end{tabular}
\end{center}

We now treat the case where $f=1$, that is, $\{1\}=\SL(1,p^d)\leq G_0\leq \Gamma L(1,p^d)\cong \ZZ_{p^d-1}: \ZZ_d$. More precisely, 
if $G_0$ has a $2$-transitive action of degree $r\geq 3$, then $r$ must be an odd prime dividing $p^d-1$ such that $r-1$ divides $d$ and the order of $p$ modulo $r$ is $r-1$. In that case,  the induced $2$-transitive permutation group  is $\ZZ_{r}: \ZZ_{r-1}\cong AGL(1,r)$.

If $d=2$, then $r=3$, $M_1$ has $\mathbb{F}_p$-dimension $1$, and by condition (a) of Construction \ref{constr-regN},  $G_0$ is transitive on the set of $p+1$ subspaces of dimension 1, so $p+1=|\calM|=r=3$ and  $p=2$. Since $\ZZ_{3}: \ZZ_2\cong \Gamma L(1,4)= \SL(2,2)$ this case has already been considered in the case $f\geq 2$ above, so we can now assume that $d>2$. Take $G_0=\langle \tau^i,\tau^j\sigma^t  \rangle$ in standard form as in Lemma \ref{GaLcriterion}.

We next consider the case $G_0\leq \Gamma L(1,2^6)\cong \ZZ_{63}: \ZZ_6$. Since $r$ is prime and divides $2^6-1=63$, $r$ must be equal to $3$ or $7$.  If $r=7$, then the order of $2$ modulo $r$ is $3\neq r-1$.  So $r=3$ and the $2$-transitive action is isomorphic to $\ZZ_{3}: \ZZ_{2}$. Thus $t=1$ or $3$. Then it follows from  Lemma \ref{GaLcriterion} that $i=1$, and so $G_0=\Gamma L(1,2^6)$ or $\langle \tau,\sigma^3\rangle\cong\ZZ_{63}: \ZZ_2$.  However in both cases the stabiliser of an element in the 2-transitive action, namely $\langle \tau^3,\sigma\rangle$ or  $\langle \tau^3,\sigma^3\rangle$ respectively, does not stabilise a subspace (it has subdegrees $1$, $21$ and $42$).

This leaves the case where $d\geq 3$ and $(p,d)\neq (2,6)$, and here, by Zsigmondy's Theorem \cite{zsi}, $p^d-1$ has a primitive prime divisor $s$, that is, $s$ is prime, $s\div p^d-1$ and $s\notdiv p^c-1$ for any $c<d$. Hence the order of $p$ modulo $s$ is $d$ and so $s-1\geq d$.
Since $G$ is $2$-transitive on $V$,  $p^d-1$ divides the order of $G_0$, and so $s\div |G_0|$. Therefore $G_0$ has elements of order $s$. Let $g$ be such an element.
Since the stabiliser of a subspace in $\GL(d,p)$ has order coprime with $s$, $g$ cannot fix any subspace of $V$.
So if $\calM=\{M_1,\ldots,M_r\}$ ($r\geq 3$) satisfies the conditions of Construction \ref{constr-regN}, then $s$ divides $r$, and since $r$ is a prime we must have $r=s$. Thus $d\leq s-1=r-1\leq d$, and so $r=s=1+d$.
Moreover  the kernel of the action of $G_0$ on $\calM$ has order coprime with $s$.
Since the induced action of $G_0$ on $\calM$ is $\ZZ_{s}: \ZZ_{s-1}$, whose order is not divisible by $s^2$,  $s^2\notdiv p^d-1$, and in fact $s$ is equal to the product, including multiplicities, of all the primitive prime divisors of $p^d-1$. 
By Hering \cite[Theorem 3.8]{Hering}  (see also \cite[Proposition 2.2 ]{NP} ), it follows that 
$(p,d)$ lies in $\{(2,4),(2,10),(2,12),(2,18),(3,4),(3,6),(5,6)\}$. 

Since the induced action of $G_0$ on $\calM$ is $\ZZ_{s}: \ZZ_{s-1}$, we have that $s-1=d$ divides $|G_0/G_0\cap\langle \tau \rangle|=|\sigma^t|=d/t$. This implies that $t=1$. Moreover $s=d+1$ divides $|G_0\cap\langle \tau \rangle|=|\langle \tau^i\rangle|=\frac{p^d-1}{i}$. In other words $i|\frac{p^d-1}{d+1}$.

In every specific case, we determine the transitive subgroups by showing that if $i>1$ then very few values of $i$ satisfy $i|\frac{p^d-1}{d+1}$ and conditions (B) and (3) of Lemma \ref{GaLcriterion}.
The following assertion  is easy to prove:
\medskip

\begin{minipage}{13cm}
 If $i=i_1i_2$ with $i_1,i_2$ coprime and $i_1|j\frac{p^{k_1}-1}{p-1}$ for some $k_1\leq i_1$ and $i_2|j\frac{p^{k_2}-1}{p-1}$ for some $k_2\leq i_2$ and at least one of the inequalities $k_1< i_1,k_2< i_2$ holds, then $i|j\frac{p^{k_1k_2}-1}{p-1}$ and $k_1k_2<i=i_1i_2$. 
\end{minipage}\hspace{1cm}
\begin{minipage}{1cm}
 (*)
\end{minipage}

\medskip
Thus if the hypotheses of (*) hold for $i$, then $i$ does not satisfy condition (B) of Lemma \ref{GaLcriterion}. 
It follows that if condition (B) fails for the $b$-part of $i$, for each prime $b$ dividing $i$, then condition (B) fails for $i$.
In particular, if no prime power $i$ satisfies condition (B) then no $i$ satisfies condition (B). In this case, we must have $i=1$ and this implies that 
$G_0=\GammaL(1,p^d)$.

Assume $(p,d)$ is one of $(2,4),(2,10),(2,12),(2,18)$, so $r$ is $5$, $11$, $13$, $19$ respectively. 
Considering for $i$ all the prime powers dividing $\frac{2^d-1}{d+1}$, which are  $3,9,27,5,7,31,73$, and taking for $k$ the values   $2,6,18,4,3,5,18$ respectively (regardless of the value of $j$), we see that condition (B) fails, so by our observation above  $i=1$ and $G_0= \Gamma L(1,2^d)$. 
A stabiliser $\langle \tau^r,\sigma\rangle$ of an element in a degree $r$ action has orbit sizes $1$, $\frac{2^d-1}{d+1}$  and $d\frac{2^d-1}{d+1}$ on $V$. Except for $d=4$, it follows that  $\langle \tau^r,\sigma\rangle$ does not stabilise a subspace.  
If $d=4$, the stabiliser $\langle \tau^5,\sigma\rangle$  has orbits $\{0\}$, $\{1,\epsilon^5,\epsilon^{10}\}=\FF_4^*$, $\FF\setminus \FF_4$, so it stabilises $\FF_4$ which is closed under addition and so is a subspace. We get in this way the $5$ hyperplanes of $V(2,4)$,  and  Line 3 of Table \ref{maintable2} holds.

Assume $(p,d)$ is one of $(3,4),(3,6)$ so $r$ is $5$, $7$ respectively. 
The prime powers dividing $\frac{3^d-1}{d+1}$ are  $2,4,8,16,13$. For $i$ one of $4,8,16,13$ we take $k=2,4,8,6$ respectively (regardless of the value of $j$),
and see that condition (B) fails. 
Since $2|j\frac{3^{2_1}-1}{3-1}$, assertion (*) implies that condition (B) fails for $i=26$, so  the only possible values for $i$ are $1$ and $2$. So  $G_0= \Gamma L(1,3^d)$ or the index $2$ subgroup $\langle \tau^2,\tau\sigma\rangle$.  A stabiliser $H$ of an element in a degree $r$ action is $\langle \tau^r,\sigma\rangle$, or $\langle \tau^{2r},\tau\sigma\rangle$ respectively. In the first case, the $H$-orbit sizes in $V$ are $1$, $\frac{3^d-1}{d+1}$ and $d\frac{3^d-1}{d+1}$; in the second case they are $1$, $\frac{3^d-1}{d+1}$,  $d\frac{3^d-1}{2(d+1)}$ and $d\frac{3^d-1}{2(d+1)}$. So again $H$ does not stabilise a subspace. 

The last case is  $(p,d)=(5,6)$. The prime powers dividing $\frac{5^6-1}{6+1}$ are  $2,4,8,3,9,31$. 
For $i$ one of $4,8,3,9,31$ we can take $k=2,2,2,6,3$ respectively (note that for $i=4,8$ condition (3) implies that $j=2,4$ respectively).
Since $2|j\frac{5^{2_1}-1}{5-1}$, assertion (*) implies that condition (B) fails for $i=2\cdot 3$, $2\cdot 9$ and $2\cdot 31$. Hence the only possible values for $i$ are $1$ and $2$. So 
  $G_0= \Gamma L(1,5^6)$ or the index $2$ subgroup $\langle \tau^2,\tau\sigma\rangle$. A stabiliser $H$ of an element in a degree $r$ action is $\langle \tau^r,\sigma\rangle$, or $\langle \tau^{2r},\tau\sigma\rangle$ respectively. In the first case, the $H$-orbit sizes in $V$ are $1$, $\frac{5^6-1}{6+1}$ and $d\frac{5^6-1}{6+1}$; in the second case they are $1$, $\frac{5^6-1}{6+1}$,  $d\frac{5^6-1}{2(6+1)}$ and $d\frac{5^d-1}{2(6+1)}$. So again $H$ does not stabilise a subspace.

\medskip

2. $\Sp(f,q)\lhd G_0$, where $f=d/e$, $q=p^e$ and $f$ is even. 
If $f=2$, then $\SL(2,q)=\Sp(2,q)\lhd G_0$ and we already treated this case. So we can assume $f\geq 4$.
If $q=2$, then $\Sp(d,2)$ ($d$ even) has $2$-transitive actions of degree $2^{d-1}+2^{d/2-1}$ and $2^{d-1}-2^{d/2-1}$ (for $d=4$, this follows from $\Sp(d,2)\cong \Sym(6)\cong P\Sigma L(2,9)$). However a stabiliser in one of these actions fixes a hyperbolic or ellipic quadric in $V=V(d,2)$ and does not stabilise any subspace. 
Looking at the list of $2$-transitive groups, there are no other (faithful or unfaithful) $2$-transitive actions for $G_0$ of degree at least $3$.

\medskip

3.   $G_2(2^e)\lhd G_0$, where $p=2$ and $d=6e$.
Suppose first that $e=1$, that is, $V=V(6,2)$ and $G_2(2)\cong P\Gamma U(3,3)\lhd G_0$. Since $G'_2(2)$ is the unique minimal normal subgroup of $G_0$, $G_0$ does not have an unfaithful $2$-transitive action of degree $\geq 3$. The only faithful $2$-transitive action of $G_0$ is of degree $28$. However it is easily checked that the stabiliser of an element in that action does not stabilise a subspace of $V(6,2)$.
Suppose now $e\geq 2$. Then $G_2(2^e)$ is the unique minimal normal subgroup of $G_0$ and $G_0/G_2(2^e)$ is cyclic, so $G_0$ does not have an unfaithful  $2$-transitive action of degree $\geq 3$. So if $G_0$ has a $2$-transitive action, it must be faithful and almost simple with socle $G_2(2^e)$, but this group has no $2$-transitive action.

\medskip

Conversely, suppose that  $\calD$ is the design $AG(f,q)$ described in Example \ref{affine} and  $G$ is as in Lines 1, 2, 3 of Table \ref{maintable2}.
By construction, $\calD$ is a non-trivial 2-design and $G$ is 2-transitive of affine type on $\calP$, in each case.  
By Lemma \ref{conv-facts}, it is sufficient to check that conditions (a), (b), (c) of Lemma \ref{facts} hold for $G$ to conclude that $\calD$ is $G$-pairwise transitive.
All these groups are 2-transitive on $\calP$, so (a) holds. Since $G$ permutes the parallel classes of hyperplanes, $G^\calB$ is imprimitive and in each case $G^\calB$ has rank 3. Also $\calD$ has non-intersecting blocks (parallel hyperplanes), and so (b) holds. Finally, since $G_0$ permutes the parallel classes of hyperplanes in each case, $G_0$ is transitive on the set $[0]$ of blocks incident with $0$ (the $(f-1)$-dimensional subspaces). Moreover, in each case, the stabiliser in $G_0$ of a hyperplane through $0$ is transitive on the other hyperplanes parallel to it, thus  $G_0$ is also transitive on the set  $\calB\setminus [0]$. We conclude that (c) holds. 
\qed

\subsubsection{ $G$  is almost simple on points}\label{sec:AS}
We study in this section case (2) of Lemma \ref{imprimcase}.
Here $G^\calB$ is imprimitive but quasiprimitive of rank 3. Again there is a unique block system $\calS$ for the $G$-action on $\calB$, and since $G^\calB$ is quasiprimitive of rank $3$, we have $G^\calS\cong G$  2-transitive.

\begin{lemma}\label{lem:noPSL}
 Suppose  $\calD$  is a non-trivial $G$-pairwise transitive 2-design, $G^\calP$ is of almost simple type, and $G^\calB$ is imprimitive of rank 3,  with block system $\calS$.
Then the two actions $G^\calP$ and $G^\calS$ are not the point and hyperplane actions of a projective linear group.\end{lemma}
\proof
Suppose the socle of $G$ is $\PSL(a,q)$, with point action on $\calP$ and hyperplane action on $\calS$ for some $a\geq 3$ and some prime power $q$. Then for $S\in \calS$, $G_S$ has two orbits on $\calP$, of sizes $\frac{q^{a-1}-1}{q-1}$ and $\frac{q^{a}-q^{a-1}}{q-1}$ (points lying in and out of the hyperplane).

Let  $S\in \calS$ and $b\in S$. Since $G_b$ also has two orbits on $\calP$  by Lemma \ref{facts}(d), and $G_b\leq G_S$, it follows that $G_S$ and $G_b$ have
 the same two orbits on $\calP$. Let $b'$ be another block in $S$. Then the same argument shows that $G_{b'}$ and $G_S$, and hence also $G_b$, have the same two orbits on $\calP$.

Since $k=|[b]|=|[b']|\leq v/2$ by Lemma \ref{lem:orbitsizes}, it follows that $k=\frac{q^{a-1}-1}{q-1}<v/2$ and $[b]=[b']$, contradicting Lemma \ref{facts}(e). 
\qed

We now give an example.
\begin{example}\label{Golay}{\rm
 Let $C$ be the ternary Golay code. Then $C$ has 24 total words (weight 12 words), denoted by $ \varepsilon w_i$ where $i\in \{0,1,2,\ldots, 10\} \cup\{\infty\}$ and  $\varepsilon\in\{1,-1\}$. The set of total words is  preserved by $M_{12}$,  and $M_{11}$ is the setwise stabiliser in $M_{12}$ of $\{w_\infty,-w_\infty\}$ where $w_\infty$ is the all-one vector (see the Atlas \cite[p.32 and 18]{atlas}). The remaining 22 total words each have  six $1$ and six $-1$ entries. 

Let $\calP=\{0,1,2,\ldots, 10\} \cup\{\infty\}$, and let $\calB=\{\varepsilon w_i|i\neq\infty, \varepsilon\in\{1,-1\}\}$, where point $j$ is incident with block $\varepsilon w_i$ exactly when the $j$-entry of $\varepsilon w_i$ is equal to 1.
 Then $\calD=(\calP,\calB,I)$ forms a 
$2-(12,6,5)$ design (also a $3-(12,6,2)$ design actually) $H(12)$. Blocks intersect in 0 or 3 points, and are partitioned into parallel classes of size 2 ($\{w_i,-w_i\}$ for each $i$),  and $\Aut(\calD)=M_{11}$.}
\end{example}
This is also the unique Hadamard 3-design associated with a $3\times 3$ Hadamard matrix characterised in \cite{Norman}. The number of points incident with two intersecting blocks, the parameters of the design, and the full automorphism group follow from Norman's work in \cite{Norman}.

\begin{theorem}\label{blockaction-rank-3-imp}
Assume  $\calD$ is a non-trivial  $2$-design and assume $G\leq \Aut({\cal D})$ is such that  $G^\calB$ has rank $3$ and  is imprimitive, and $G^\calP$ is of almost simple type.
Then $\calD$ is $G$-pairwise transitive if and only if either
 $\calD$ is $\AG(3,2)$ (see Example \ref{affine}) or $\calD$ is as in Example \ref{Golay}, and $G$ is $\PSL(2,7)$ or $M_{11}$ as in Line 4 or 5 of Table \ref{maintable2}, respectively.

\end{theorem}

\begin{center}
\begin{table}[ht]
\begin{tabular}{|l|l|l|l|p{7cm}|}
\hline
$G$ & $|\calS|$ & $|S|$&  $G_S^S$& extra conditions \\
\hline
$M_{11}$& $11$&$2$&$C_2$&\\
$G\geq \PSL(2,q)$&$q+1$&$2$&$C_2$& $q\equiv 1 \pmod 4$, $q$ not a prime and $G$  satisfies the conditions explained in \cite[Proposition 5.10]{rk3} \\
$G\geq \PSL(a,q)$& $\frac{q^a-1}{q-1}$&$m$& $\AGL(1,m)$&  $a\geq 3$, $q\geq 3$ and $G$  satisfies the conditions described in  \cite[Table 1]{rk3}\\
$\PGL(3,4)$&$21$&$6$&$\PSL(2,5)$&\\
$\PGammaL(3,4)$&$21$&$6$&$\PGL(2,5)$&\\
$\PSL(3,5)$&$31$&$5$& $\Sym(5)$&\\
$\PSL(5,2)$& $31$&$8$&$\Alt(8)$&\\
$\PGammaL(3,8)$&$73$&$28$&  $\Ree(3)$&\\
$\PSL(3,2)$&$7$&$2$& $C_2$&\\
$\PSL(3,3)$&$13$&$3$&$\Sym(3)$&\\
 \hline
\end{tabular}
\caption{Quasiprimitive imprimitive rank 3 groups} \label{bigtable3}
\end{table}
\end{center}

\vspace*{-1cm}
\proof
Assume  $\calD$ is  $G$-pairwise transitive. Let $\calS$ be the unique nontrivial system of imprimitivity of $G$ on blocks (see \cite[Lemma 2.3]{rk3}). Since $G^\calB$ is quasiprimitive by Lemma \ref{imprimcase}, it follows that $G\cong G^\calS$ by  \cite[Corollary 2.5]{rk3}. 
A list of imprimitive but quasiprimitive rank 3 groups $G^\calS$, together with $|\calS|$ and $|S|$, $G_S^S$ and  conditions, is given in \cite[Table 1]{rk3} (see \cite[Theorem 1.2]{rk3}) and reproduced here for convenience in Table \ref{bigtable3} (at the end of the paper), where we write more specifically the conditions we are going to need in our proof. In particular, for the third row, the conditions (see Line 12 of \cite[Table 3]{rk3}) imply that $q$ cannot be equal to $2$, since $m$ must be prime and $md$ is a divisor of $q-1$.

By Lemmas \ref{equiv2trans} and \ref{lem:noPSL}, $G$ must have two non-equivalent 2-transitive actions (on $\calP$ and on $\calS$), and if $G$ has socle $\PSL(a,q)$ ($a\geq 3$) then these are not the actions on points and hyperplanes of the projective geometry. 
The groups $G$ in  Table \ref{bigtable3} which satisfy these additional properties are very few:
$M_{11}$  of degree 11 on $\calS$ and degree 12 on $\calP$; $\PSL(2,9)$ and $\PSigmaL(2,9)$ of degree 10 on $\calS$ and degree 6 on $\calP$; and $\PSL(3,2)$ of degree 7 on $\calS$ and degree 8 on $\calP$  (information from the  Atlas \cite{atlas} was used to determine the overgroups for $\PSL(2,9)$ and $\PSL(3,2)$).


%

Suppose  $G\geq \PSL(2,9)$, with natural 2-transitive action of degree $10$ on $\calS$ and  2-transitive action of degree $6$ on $\calP$.  
The conditions \cite[Proposition 5.10]{rk3}  imply that $G=\langle \PSL(2,9),\sigma\tau\rangle\cong M_{10}$, where $\sigma$ is the Frobenius automorphism and 
$\tau=
\begin{pmatrix}
 1&0\\
0& \omega
\end{pmatrix}$
  ($\omega$ is a primitive element of $\GF(9)$). This is a contradiction because $M_{10}$ has no 2-transitive degree 6 action (see \cite[p.4]{atlas}). So none of these groups satisfy all the required conditions.



Suppose  $G= \PSL(3,2)$, with natural 2-transitive action of degree $7$ on $\calS$ and  2-transitive action of degree $8$ on $\calP$. 
For $b\in \calB$, $G_b$ has two orbits of size 4 on $\calP$. We choose one of these orbits as the set $[b]$ of points incident with $b$ and we get a $2-(8,4,3)$ design $\calD$ (this follows from an easy counting). Note that if we had chosen the other orbit, we would have obtained a design isomorphic to $\calD$ (switching the incidence between every two blocks in the same part of $\calS$ yields an isomorphism). There are four designs with these parameters up to isomorphism, but only one of them has disjoint blocks \cite[p.27]{HCD}:  namely $\AG(3,2)$ (Example \ref{affine} for $(f,q)=(3,2)$), thus $\calD=\AG(3,2)$ and Line 4 of Table \ref{maintable2} holds.

Suppose  $G=M_{11}$, with natural 2-transitive action of degree $11$ on $\calS$ and 2-transitive action  of degree $12$ on $\calP$. 
For $b\in \calB$, $G_b$ has two orbits of size 6 on $\calP$. Choosing either of these orbits  as the set $[b]$ gives by an easy counting a $2-(12,6,5)$ design $\calD$. Since $M_{11}$ is actually 3-transitive of degree 12, $\calD$ is also a 3-design with parameters $3-(12,6,2)$. 
Its derived design is a $2-(11,5,2)$ design, which is a Hadamard 2-design, shown to be unique by Todd \cite[Section 3]{Todd33}. It follows by \cite[Lemma 1]{Norman} (see also \cite[Corollary II.8.11]{BJL}) that there is a unique  
 $3-(12,6,2)$ design:  namely $H(12)$, described in Example \ref{Golay}.  Thus $\calD=H(12)$ and  Line 5 of Table \ref{maintable2} holds.

We now prove the converse.
Let $\calD$ and $G$ be as in Line 4 or 5 of Table \ref{maintable2}. Note that in both cases, the complement of a block is also a block. Moreover the set of complementary block pairs  forms a system of imprimitivity for the action on $\calB$, so  $G^\calB$ is imprimitive.
By construction, $G^\calP$ has rank 2 (so condition (a) of Lemma \ref{facts} holds), $G$ is simple, and $G^\calB$ is quasiprimitive (since  $G\cong G^\calB$ is simple and transitive).  Moreover $G^\calB$ has rank 3 by \cite{rk3}.
We checked for each case  by computer that a block stabiliser in $G$ has two orbits on  $\calP$: the points in the block, and  the points outside the block.
 Hence condition (d) of Lemma \ref{facts} holds. 
To verify condition (b) of Lemma \ref{facts}, we need to check that there are non-intersecting blocks. For $\calD=AG(3,2)$, this is because the design has parallel planes. For $\calD=H(12)$, it is true by construction, and is also proved in \cite[Lemma 1]{Norman}. Thus (b) holds.
Hence by Lemma \ref{conv-facts}, $\calD$ is $G$-pairwise transitive. This completes the proof
\qed

The proof of Theorem \ref{main} follows from  Proposition \ref{symmetric2d}, and Theorems \ref{blockaction-rk3-prim}, \ref{rank3impaffine}, \ref{blockaction-rank-3-imp}.

\section{Proof of Corollary \ref{cor}}

Suppose $\calD$ is a non-trivial  $2$-$(v,k,\lambda)$ design, and $G\leq \Aut({\cal D})$ such that  $\calD$ is $G$-pairwise 
transitive and $N$-nicely affine for some non-trivial normal subgroup $N$ of $G$.
By \cite[Theorem 1.5]{starquo}, $G^\calB$ is imprimitive and of rank 3. 
Then by Lemma \ref{imprimcase} $G$ is either of affine type or of almost simple type. If $G$ is of almost simple type, then $G^\calB$ is quasiprimitive, so  all normal subgroups $N$ of $G$ are transitive on $\calB$, and hence the orbits of $N$ cannot be parallel classes. 
Thus $G$ is of affine type and $\calD$, $G$ are as in one the first three Lines of Table \ref{maintable2}.

Conversely, assume  $\calD$, $G$ are as  in one the first three lines of Table \ref{maintable2}. 
These all satisfy $G^\calB$ being imprimitive and of rank 3, and with $G$ 2-transitive of affine type on $\calP$, so by  Lemma \ref{imprimcase}(1), they are $N$-nicely affine for  $N$ the translation subgroup of $G$.

\end{document}